\documentclass[english,reqno]{amsart}
\usepackage{times}
\usepackage[T1]{fontenc}
\usepackage{setspace}
\usepackage{amssymb}

\makeatletter
 \theoremstyle{plain}    
 \newtheorem{thm}{Theorem}[section]
 \numberwithin{equation}{section} 
 \numberwithin{figure}{section} 
 \theoremstyle{plain}
 \theoremstyle{remark}
 \newtheorem*{rem*}{Remark}

\usepackage{psfrag}
 \usepackage{bbm}

\renewcommand{\tilde}{\widetilde}

\newcommand{\Z}{\mathbb{Z}}

\newcommand{\Q}{\mathbb {Q}}

\renewcommand{\H}{\mathbb {H}}

\newcommand{\1}{\boldsymbol{1}}

\newcommand{\N}{\mathbb {N}}

\newcommand{\I}{\mathbb {I}}

\AtBeginDocument{
  
}

\usepackage{babel}
\makeatother
\begin{document}
\title{Large deviation asymptotics for continued fraction expansions}

\author{Marc Kesseböhmer and Mehdi Slassi}

\address{Fachbereich 3 - Mathematik und Informatik, Universit\"{a}t Bremen,
D--28359 Bremen, Germany}

\email{mhk@math.uni-bremen.de, slassi@math.uni-bremen.de}

\begin{abstract}
We study large deviation asymptotics for processes defined in terms
of continued fraction digits. We use the continued fraction digit
sum process to define a stopping time and derive a joint large deviation
asymptotic for the upper and lower fluctuation process. Also a large
deviation asymptotic for single digits is given.
\end{abstract}



\thanks{The research was partly supported by the ZF University of Bremen
(Grant No. 03/106/2).}

\maketitle

\section{Introduction }

\let\languagename\relax Any number $x\in\I:=\left[0,1\right]\setminus\Q$
has a simple infinite continued fraction expansion\[
x=\frac{1}{\kappa_{1}\left(x\right)+{\displaystyle \frac{1}{\kappa_{2}\left(x\right)+\cdots}}},\]
 where the unique \emph{continued fraction digits} $\kappa_{n}\left(x\right)$
are from the positive integers $\mathbb{N}$. This expansion is closely
related to the \emph{Gauss transformation $G:\I\rightarrow\I$} given
by\[
G(x):=\frac{1}{x}-\left\lfloor \frac{1}{x}\right\rfloor ,\]
 where $\left\lfloor x\right\rfloor $ denotes the greatest integer
not exceeding $x\in\mathbb{R}$. We write $G^{n}$ for the $n$-th
iterate of $G$, $n\in\mathbb{N}_{0}=\left\{ 0,1,2,\ldots\right\} $,
with $G^{0}:=\textrm{id.}$ Then for all $n\in\mathbb{N}$ we have
$\kappa_{n}(x)=\left\lfloor \left(G^{n-1}\left(x\right)\right)^{-1}\right\rfloor .$
Clearly, $\left(\kappa_{n}\right)$ defines a sequence of random variables
on the measure space $\left(\I,\mathcal{B},\mathbb{P}\right)$, where
$\mathcal{B}$ denotes the Borel $\sigma$-algebra of $\I$ and $\mathbb{P}$
some probability measure on $\mathcal{B}$. With respect to $\lambda$,
denoting the Lebesgue measure restricted to $[0,1]$, each $\kappa_{n}$
has infinite expectation. By the ergodicity of the Gauss transformation
with respect to the Gauss measure $d\mu\left(x\right):=\frac{1}{\log2}\frac{1}{1+x}d\lambda\left(x\right)$
we readily reproduce Khinchin's result on the geometric mean of the
continued fraction digits, i.e.\[
\sqrt[n]{\kappa_{1}\cdots\kappa_{n}}\rightarrow K,\quad\lambda\textrm{-a.e.,}\]
 where $K=2,685...$ denotes the Khinchin constant. This in particular
implies that \[
\liminf_{n}\frac{\kappa_{n}}{n}=0\quad\lambda\textrm{-a.e.}\]
 Also by a classical result of Khinchin (cf. \cite{Khinchin:36}),
we know that for $\lambda$-almost every $x\in\left[0,1\right]$ we
have for infinitely many $n\in\N$ that $\kappa_{n}\left(x\right)>n\log n.$
Let $S_{n}(x):=\kappa_{1}(x)+\cdots+\kappa_{n}\left(x\right)$, $x\in\I$,
denote the sum of the first $n$ digits. Then Khinchin result implies
in particular that $\lambda$--a.e. \[
\lim_{n\rightarrow\infty}\frac{S_{n}}{n}=+\infty\quad\textrm{ and }\quad\limsup_{n\to\infty}\frac{\kappa_{n}}{n}=+\infty.\]
 Diamond and Vaaler have shown in \cite{DiamondVaaler:86} that for
the trimmed sum $S_{n}^{\flat}:=S_{n}-{\displaystyle \max_{1\leq\ell\leq n}\kappa_{\ell}}$
we have\[
\lim_{n\rightarrow\infty}\frac{S_{n}^{\flat}}{n\log n}=\frac{1}{\log2}\;\;\lambda\textrm{-a.e.}\]
 This demonstrates that the intricate stochastic properties of $S_{n}$
arise from the occurrences of rare but exceptionally large continued
fraction digits. However, in order to find a suitable stochastic description
of the processes $S_{n}$ and $\kappa_{n}$ we define for $n\in\N$
the following (stopping time) process\[
\theta_{n}\left(x\right):=\max\left\{ \ell\in\mathbb{N}_{0}\;:\;\sum_{i=1}^{\ell}\kappa_{i}\left(x\right)\leq n\right\} ,\quad x\in\I.\]
Employing infinite ergodic theory it has been shown in \cite{KessboehmerSlassi:05}
that with respect to $\lambda$ the following convergence in distribution
holds\[
\frac{\log\left(n-S_{\theta_{n}}\right)}{\log(n)}\stackrel{\lambda}{\;\Longrightarrow\boldsymbol{U}},\]
where $\boldsymbol{U}$ denotes the random variable uniformly distributed
on the unit interval. Furthermore, $n-S_{\theta_{n}}$ obeys a large
deviation asymptotic, i.e. for all $x\in\left(0,1\right)$ we have\begin{equation}
\lambda\left(\frac{n-S_{\theta_{n}}}{n}>x\right)\sim\frac{-\log\left(x\right)}{\log\left(n\right)}\qquad\textrm{as}\: n\to\infty.\label{eq:LDold}\end{equation}
In here, $a_{n}\sim b_{n}$ stands for $\lim\frac{a_{n}}{b_{n}}=1$. 

The following list gathers the new large deviation laws valid for
the continued fraction expansion. Note that the result in (1) generalizes
the asymptotic in (\ref{eq:LDold}).

\begin{enumerate}
\item For $0\leq x<1$ and $y\geq0$ with $x+y\neq0$ we have\[
\lambda\left(\frac{n-S_{\theta_{n}}}{n}>x,\frac{S_{\theta_{n}+1}-n}{n}>y\right)\sim\log\left(\frac{1+y}{x+y}\right)\cdot\frac{1}{\log\left(n\right)}\qquad\textrm{as}\: n\to\infty.\]
 
\item For $x>0$ we have \begin{equation}
\lambda\left(\frac{\kappa_{\theta_{n}+1}}{n}>x\right)\sim\frac{H\left(x\right)}{\log\left(n\right)}\quad\textrm{as}\quad n\to\infty,\label{gl2}\end{equation}
where $H$ is the convex, decreasing, differentiable function given
by\[
H\left(x\right):=\left\{ \begin{array}{ll}
1-\log\left(x\right) & \quad\textrm{for }\: x\in\left(0,1\right),\\
1/x & \quad\textrm{for }\: x\geq1.\end{array}\right.\]

\item For $x\in\left(0,1\right)$ we have \begin{eqnarray*}
\lambda\left(\frac{\kappa_{\theta_{n}+1}}{\sum_{k=1}^{\theta_{n}+1}\kappa_{k}}>x\right) & \sim & \frac{\tilde{H}\left(x\right)}{\log\left(n\right)}.\end{eqnarray*}
where $\tilde{H}$ is a decreasing, differentiable function given
by\[
\tilde{H}\left(x\right):=\frac{1-x}{x}\log\left(\frac{1}{1-x}\right)+\log\left(\frac{1}{x}\right)\quad\textrm{for }\: x\in\left(0,1\right).\]

\item For $x>0$ we have \[
\lambda\left(\frac{\kappa_{\theta_{n}+1}}{\sum_{k=1}^{\theta_{n}}\kappa_{k}}>x;\theta_{n}>0\right)\sim\frac{\tilde{H}\left(\frac{x}{1+x}\right)}{\log\left(n\right)}.\]

\end{enumerate}
\begin{rem*}
We would like to remark that using for example results from \cite{KesseboehmerStratmann:04b,KesseboehmerStratmann:07}
it is possible to give an interpretation of the above processes in
terms of geodesic windings on the cusped orbifold $\H/\textrm{PSL}_{2}\left(\Z\right)$
with respect to visits to its compact part. This shows that also smooth
dynamical systems can be investigated by the renewal theoretical methods
developed in this paper. 

As a corollary of the results in \cite{KesseboehmerSlassi:07} we
also get\[
\frac{\log\left(S_{\theta_{n}+1}-n\right)}{\log\left(n\right)}\stackrel{\lambda}{\;\Longrightarrow}\;\boldsymbol{U},\quad\textrm{and}\quad\frac{\log\left(\kappa_{\theta_{n}+1}\right)}{\log\left(n\right)}\stackrel{\lambda}{\;\Longrightarrow}\;\boldsymbol{U}.\]
This will be explained at the end of Subsection \ref{sub:Renewal-theory-for}.
\end{rem*}
Related results on the digit sum process can also be found in \cite{Heinrich:87},
\cite{Hendley:00}, and \cite{GuivarchLeJan:93,GuevarchLeJan:96}.

The paper is organized as follows. In Section \ref{sec:Preliminaries}
we first introduce the necessary concepts from infinite ergodic theory.
Then, in Subsection \ref{sub:Renewal-theory-for} we link the number
theoretical processes under consideration to the renewal theory defined
with respect to the Farey map for which the general results from infinite
ergodic theory are applicable . 

In Section \ref{sec:Large-Deviation} we state the main theorem in
its general form and then give its proof.

\section{Preliminaries\label{sec:Preliminaries}}

\subsection{Infinite ergodic theory. \label{sub:Infinite-ergodic-theory} }

In this subsection we briefly recall the basic concepts and result
from infinite ergodic theory as needed in this paper. For the definition
and further details we refer the reader to \cite{Aaronson:97}. Let
$\left(X,T,\mathcal{A},\mu\right)$ be a conservative ergodic measure
preserving dynamical systems where $\mu$ is an infinite $\sigma$-finite
measure. Let \[
\mathcal{P}_{\mu}:=\left\{ \nu:\nu\:\textrm{probability measure on}\:\mathcal{A}\,\textrm{with }\nu\ll\mu\right\} \]
 denote the set of probability measures on $\mathcal{A}$ which are
absolutely continuous with respect to $\mu$. The measures from $\mathcal{P}_{\mu}$
represent the admissible initial distributions for the processes under
consideration. With $\mathcal{P}_{\mu}$ we will sometimes also denote
the set of the corresponding densities. 

Let us recall the notion of the wandering rate. For a fixed set $A\in\mathcal{A}$
with $0<\mu\left(A\right)<\infty$ we set for all $n\geq0$

\[
A_{n}:=\bigcup_{k=0}^{n}T^{-k}A\quad\mathrm{and}\quad W_{n}:=W_{n}\left(A\right):=\mu\left(A_{n}\right)=\sum_{k=0}^{n}\mu\left(A\cap\{\varphi>k\}\right),\]
and call the sequence $\left(W_{n}\left(A\right)\right)$ the \emph{wandering
rate} of $A.$ In here, \begin{equation}
\varphi(x)=\inf\{ n\geq1:\; T^{n}(x)\in A\},\;\; x\in X,\label{returntime}\end{equation}
 denotes the \emph{first return time} to the set $A$. Since $T$
is conservative and ergodic, for all $\nu\in\mathcal{P_{\mu}}$, we
have $\nu\left(\left\{ \varphi<\infty\right\} \right)=1.$

To explore the stochastic properties of $T$ it is often useful to
study the long-term behaviour of the iterates of its \emph{transfer
operator} $\hat{T}$, i.~e. $\hat{T}:L_{1}\left(\mu\right)\rightarrow L_{1}\left(\mu\right)$
is the dual of \emph{}$f\mapsto f\circ T$, $f\in L_{\infty}\left(\mu\right)$. 

The ergodic properties of $(X,T,\mathcal{A},\mu)$ can be characterized
in terms of the transfer operator in the following way. $T$ is conservative
and ergodic if and only if for all $f\in L_{1}^{+}\left(\mu\right):=\left\{ f\in L_{1}\left(\mu\right):\; f\geq0\;\mathrm{and}\;\int_{X}f\; d\mu>0\right\} $
we have $\mu$-a.e. $\sum_{n\geq0}\hat{T}^{n}\left(f\right)=\infty$. 

A set $A\in\mathcal{A}$ with $0<\mu\left(A\right)<\infty$ is called
\emph{uniform for} $f\in\mathcal{P_{\mu}}$ if there exists a sequence
$\left(a_{n}\right)$ of positive reals such that\[
\frac{1}{a_{n}}\sum_{k=0}^{n-1}\hat{T}^{k}\left(f\right)\;\longrightarrow\;1\qquad\mu-\textrm{a.e. \; uniformly\; on}\; A\]
The set $A$ is called a \emph{uniform} set \emph{}if it is uniform
for some $f\in\mathcal{P_{\mu}}.$

Note that from \cite{Aaronson:86} we know, that $\left(a_{n}\right)$
is regularly varying with exponent $\alpha$ (for the definition of
this property see \cite{BinghamGoldieTeugels:89}) if and only if
$\left(W_{n}\right)$ is regularly varying with exponent $\left(1-\alpha\right)$.
In this case $\alpha$ lies in the interval $\left[0,1\right]$ and
\begin{equation}
a_{n}W_{n}\sim\frac{n}{\Gamma\left(1+\alpha\right)\Gamma\left(2-\alpha\right)}.\label{eq:AsymptAaronson}\end{equation}

Next, we recall the notion of uniformly returning sets, which will
be used to proof the statements in Theorem \ref{thm:neue2} (cf. \cite{KesseboehmerSlassi:07})

\begin{itemize}
\item A set $A\in\mathcal{A}$ with $0<\mu\left(A\right)<\infty$ is called
\emph{uniformly returning} \emph{for} $f\in\mathcal{P_{\mu}}$ if
there exists an positive increasing sequence $\left(b_{n}\right)$
diverging to $\infty$ such that\[
b_{n}\hat{T}^{n}\left(f\right)\;\longrightarrow\;1\quad\mu-\textrm{a.e. \; uniformly\; on}\; A.\]

\item The set $A$ is \emph{}called \emph{uniformly} \emph{returning} if
it is uniformly returning for some $f\in\mathcal{P_{\mu}}.$
\end{itemize}
In \cite{KesseboehmerSlassi:07} it is shown that every uniformly
returning set $A$ is uniform but not necessary with respect to the
same function. If $A$ is uniformly returning for $f\in\mathcal{P}_{\mu}$
such that $f$ is bounded, then $A$ is also uniform for $f$. This
observation will be relevant in Theorem \ref{thm:neue2}. 

From \cite{KesseboehmerSlassi:07}, we know that $\left(b_{n}\right)$
is regularly varying with exponent $\beta\in\left[0,1\right)$ if
and only if $\left(W_{n}\right)$ is regularly varying with the same
exponent. In this case,\[
b_{n}\sim W_{n}\Gamma\left(1-\beta\right)\Gamma\left(1+\beta\right)\qquad\left(n\to\infty\right).\]
The following fact will be crucial in the proof of the Theorem \ref{thm:neue2}.

\begin{itemize}
\item [(UA)] \textbf{Uniform asymptotic} (\cite{Seneta:76}) Let $\left(p_{n}\right)$
and $\left(q_{n}\right)$ be two positive sequences with $p_{n}\rightarrow\infty$
and $\frac{p_{n}}{q_{n}}\in\left[1/K,K\right],\; K\geq1$ for $n$
large enough. Then for every slowly varying function $L$ we have\[
\lim_{n\to\infty}\frac{L\left(p_{n}\right)}{L\left(q_{n}\right)}=1.\]

\end{itemize}

\subsection{Renewal theory for continued fractions\label{sub:Renewal-theory-for}}

In this subsection we connect the number theoretical processes considered
in this paper with the renewal processes defined with respect to the
Farey map.

Let $T:\left[0,1\right]\rightarrow\left[0,1\right]$ be the \emph{Farey
map}, defined by\[
T\left(x\right):=\left\{ \begin{array}{ll}
T_{0}\left(x\right), & x\in\left[0,\frac{1}{2}\right],\\
T_{1}\left(x\right), & x\in\left(\frac{1}{2},1\right],\end{array}\right.\]
where\[
T_{0}\left(x\right):=\frac{x}{1-x}\qquad\textrm{and}\qquad T_{1}\left(x\right):=\frac{1}{x}-1.\]
Then $\left(\left[0,1\right],T,\mathcal{B},\mu\right)$ defines a
conservative ergodic measure preserving dynamical system, where $\mu$
denotes the $\sigma$--finite invariant measure with density $h\left(x\right):=\frac{d\mu}{d\lambda}\left(x\right)=\frac{1}{x}$
(see \cite{Thaler:83}).

Let $K_{1}:=\left(\frac{1}{2},1\right].$ Then for $n\geq0$ we consider
the following processes given by 

\begin{itemize}
\item $Z_{n}(x):=\left\{ \begin{array}{ll}
\max\left\{ k\leq n:\; T^{k}(x)\in K_{1}\right\} , & x\in A_{n}:=\bigcup_{k=0}^{n}T^{-k}K_{1},\\
0, & \textrm{else.}\end{array}\right.$
\item $Y_{n}\left(x\right):=\min\left\{ k>n:\; T^{k}\left(x\right)\in K_{1}\right\} ,\quad x\in X,$
\item $V_{n}:=Y_{n}-Z_{n}.$
\end{itemize}
In terms of renewal theory the process $n-Z_{n}$ is called \emph{spent
waiting time} since the last visit to $A$, $Y_{n}-n$ the \emph{residual
waiting time} to the next visit to $A$, and $V_{n}$ the \emph{total
waiting time} between to visits to $A$. As $T$ is conservative and
ergodic, we have for all $\nu\in\mathcal{P_{\mu}}$\[
\lim_{n\to\infty}\nu\left(A_{n}\right)=1,\,\quad\textrm{and}\quad\nu\left(\left\{ Y_{n}<\infty\right\} \right)=1\;\textrm{for\; all}\; n\geq1.\]
Furthermore, for $1\leq k\leq n\leq m$ we have (see \cite{Dynkin:61})\begin{equation}
\left\{ Z_{n}\leq k,\; Y_{n}>m\right\} =\left\{ Z_{m}\leq k\right\} .\label{Dinkin}\end{equation}
From \cite{KessboehmerSlassi:05} we know that the set $K_{1}$ is
uniformly returning for any function $f$ from the set\[
\mathcal{D}:=\left\{ f\in\mathcal{P}_{\mu}:\; f\in\mathcal{C}^{2}\left(\left(0,1\right)\right)\;\textrm{with}\; f'>0\;\textrm{and}\; f''\leq0\right\} .\]
For the wandering rate we have \[
W_{n}:=W_{n}\left(K_{1}\right)=\int_{\frac{1}{n+2}}^{1}\frac{1}{x}\, dx=\log\left(n+2\right)\sim\log\left(n\right)\qquad\left(n\to\infty\right).\]
Additional, since $K_{1}\cap\left\{ \varphi>n\right\} =\left[\frac{n+2}{n+3},1\right]$
we have\begin{equation}
\mu\left(K_{1}\cap\left\{ \varphi>n\right\} \right)=\int_{\frac{n+2}{n+3}}^{1}\;\frac{1}{x}\; dx\;\sim\;\frac{1}{n}\quad\textrm{as}\; n\to\infty.\label{Mu1}\end{equation}
 The inverse branches of the Farey map are \begin{eqnarray*}
u_{0}\left(x\right) & := & \left(T_{0}\right)^{-1}\left(x\right)={\displaystyle \frac{x}{1+x}},\\
u_{1}\left(x\right) & := & \left(T_{1}\right)^{-1}\left(x\right)={\displaystyle \frac{1}{1+x}}.\end{eqnarray*}
For $x\neq0$ the map $u_{0}\left(x\right)$ is conjugated to the
right translation $x\mapsto F\left(x\right):=x+1,$ i.e.\[
u_{0}=J\circ F\circ J\qquad\textrm{with}\quad J\left(x\right)=J^{-1}\left(x\right)=\frac{1}{x}.\]
This shows that for the $n$-th iterate we have\begin{equation}
u_{0}^{n}\left(x\right)=J\circ F^{n}\circ J\left(x\right)=\frac{x}{1+nx}.\label{u0}\end{equation}
Moreover, we have $u_{1}\left(x\right)=J\circ F\left(x\right).$ 

Let $\left\{ K_{n}:=\left(\frac{1}{n+1},\frac{1}{n}\right]\right\} _{n\geq1}$
be the countable collection of pairwise disjoint subintervals of $\left[0,1\right]$.
Setting $K_{0}=\left[0,1\right)$ it is easy to check that $T\left(K_{n}\right)=K_{n-1}$
for all $n\geq1.$ The \emph{first entry time} $e:\I\rightarrow\mathbb{N}$
into the interval $K_{1}$ is defined as\[
e\left(x\right):=\min\left\{ k\geq0:\; T^{k}\left(x\right)\in K_{1}\right\} .\]
Then the first entry and the first return time is connected to the
first digit in the continued fraction expansion by\[
\kappa_{1}\left(x\right)=1+e\left(x\right)\quad\textrm{and}\quad\varphi\left(x\right)=\kappa_{1}\circ T\left(x\right),\quad x\in\I.\]
We now consider the \emph{induced map} $S:\I\rightarrow\I$ defined
by\[
S\left(x\right):=T^{e\left(x\right)+1}\left(x\right).\]
Since for all $n\ge1$ $\left\{ x\in\I:\; e\left(x\right)=n-1\right\} =K_{n}\cap\I,$
we have by (\ref{u0}) \[
S\left(x\right)=T^{n}\left(x\right)=T_{1}\circ T_{0}^{n-1}\left(x\right)=\frac{1}{x}-n=\frac{1}{x}-\kappa_{1}(x),\quad x\in K_{n}\cap\I.\]
This implies that the induced transformation $S$ coincides with Gauss
map $G$ on $\I$.

Let $\left(\tau_{n}\right)_{n\in\mathbb{N}}$ be the sequence of return
times, i.~e. integer valued positive random variables defined recursively
by \begin{eqnarray*}
\tau_{1}\left(x\right) & := & \varphi(x)=\inf\left\{ p\geq1:\; T^{p}(x)\in K_{1}\right\} ,\quad x\in\left[0,1\right],\\
\tau_{n}\left(x\right) & := & \inf\left\{ p\geq1:\; T^{p+\sum_{k=1}^{n-1}\tau_{k}\left(x\right)}(x)\in K_{1}\right\} ,\quad x\in\left[0,1\right].\end{eqnarray*}
Then we consider the following \emph{}process, given by\[
N_{n}(x):=\left\{ \begin{array}{ll}
\max\left\{ l\leq n\;:\;\sum_{i=1}^{l}\tau_{i}\leq n\right\} , & x\in A_{n}=\bigcup_{k=0}^{n}T^{-k}K_{1},\\
0, & \textrm{else.}\end{array}\right.\]
Note that $Z_{n}=\sum_{i=1}^{N_{n}}\tau_{i}.$ Using a similarly argument
as in the proof of Lemma 3.1 in \cite{KessboehmerSlassi:05} we obtain
for all $x\in\I$ and $n\geq1$\[
\theta_{n}\left(x\right)=\left\{ \begin{array}{ll}
N_{n-1}\left(x\right)+1,\qquad & x\in K_{1},\\
N_{n-1}\left(x\right), & \textrm{else.}\end{array}\right.\]
This gives that for all $x\in\I$ and $n\geq1$\begin{equation}
S_{\theta_{n}\left(x\right)}=\left\{ \begin{array}{ll}
Z_{n-1}\left(x\right)+1,\quad & x\in A_{n-1},\\
0, & \textrm{else,}\end{array}\right.,\quad S_{\theta_{n}\left(x\right)+1}=1+Y_{n-1}\left(x\right),\label{teta1}\end{equation}
and\begin{equation}
\kappa_{\theta_{n}\left(x\right)+1}\left(x\right)=\left\{ \begin{array}{ll}
V_{n-1}\left(x\right),\qquad & x\in A_{n-1},\\
1+Y_{n-1}\left(x\right), & \textrm{else.}\end{array}\right.\label{teta2}\end{equation}
From this und the fact that $W_{n}\sim\log\left(n\right)$ we deduce
from Theorem 3.2 in \cite{KesseboehmerSlassi:07} that for all $\nu\in\mathcal{P_{\mu}}$
the following distributional convergence holds\begin{equation}
\frac{\log\left(S_{\theta_{n}+1}-n\right)}{\log\left(n\right)}\quad\textrm{and}\quad\frac{\log\left(\kappa_{\theta_{n}+1}\right)}{\log\left(n\right)}\stackrel{\mathcal{\nu}}{\;\Longrightarrow}\;\boldsymbol{U}.\label{eq:UniformLaw}\end{equation}

\section{Large Deviation\label{sec:Large-Deviation}}

We are now in the position to state our main result which are slightly
more general than the corresponding statements in the introduction.

\begin{thm}
\label{thm:neue2} For $f\in\mathcal{D}$ set $d\nu:=f\; d\mu.$ Then 
\begin{enumerate}
\item For $0\leq x<1$ and $y\geq0$ with $x+y\neq0$ we have\begin{equation}
\nu\left(\frac{n-S_{\theta_{n}}}{n}>x,\;\frac{S_{\theta_{n}+1}-n}{n}>y\right)\sim\log\left(\frac{1+y}{x+y}\right)\frac{1}{\log(n)}\qquad\textrm{as}\quad n\to\infty.\label{eqq1}\end{equation}

\item For $x>0$ we have\begin{equation}
\nu\left(\frac{\kappa_{\theta_{n}+1}}{n}>x\right)\sim\frac{H\left(x\right)}{\log(n)}\qquad\textrm{as}\quad n\to\infty,\label{eqq2}\end{equation}
where \[
H\left(x\right):=\left\{ \begin{array}{ll}
1-\log\left(x\right) & \quad\textrm{for }\; x\in\left(0,1\right),\\
1/x & \quad\textrm{for }\; x\geq1.\end{array}\right.\]

\item For $x\in\left(0,1\right)$ we have \begin{eqnarray}
\nu\left(\frac{\kappa_{\theta_{n}+1}}{\sum_{k=1}^{\theta_{n}+1}\kappa_{k}}>x\right) & \sim & \frac{\tilde{H}\left(x\right)}{\log\left(n\right)}\qquad\textrm{as}\quad n\to\infty.\label{eqq3}\end{eqnarray}
where \[
\tilde{H}\left(x\right):=\frac{1-x}{x}\log\left(\frac{1}{1-x}\right)+\log\left(\frac{1}{x}\right)\quad\textrm{for }\; x\in\left(0,1\right).\]

\item For $x>0$ we have \begin{equation}
\lambda\left(\frac{\kappa_{\theta_{n}+1}}{\sum_{k=1}^{\theta_{n}}\kappa_{k}}>x;\;\theta_{n}>0\right)\sim\frac{\tilde{H}\left(\frac{x}{1+x}\right)}{\log\left(n\right)}\qquad\textrm{as}\quad n\to\infty.\label{eqq4}\end{equation}

\end{enumerate}
\end{thm}
\begin{proof} \textbf{ad (1)} Let $0\leq x<1$ and $y\geq0$ be fixed
with $x+y\neq0$. First note, that from (\ref{teta1}) it follows
that\[
\nu\left(\frac{n-S_{\theta_{n}}}{n}>x,\;\frac{S_{\theta_{n}+1}-n}{n}>y\right)\sim\nu\left(\frac{n-Z_{n}}{n}>x,\frac{Y_{n}-n}{n}>y\right)\quad\textrm{as}\quad n\to\infty.\]
Therefore, to prove (\ref{eqq1}) it suffices to show\[
\nu\left(\frac{n-Z_{n}}{n}>x,\frac{Y_{n}-n}{n}>y\right)\sim\log\left(\frac{1+y}{x+y}\right)\frac{1}{\log(n)}\qquad\textrm{as}\quad n\to\infty.\]
In fact, by (\ref{Dinkin}) we have \begin{eqnarray*}
\nu\left(\frac{n-Z_{n}}{n}>x,\frac{Y_{n}-n}{n}>y\right) & = & \nu\left(Z_{n}\leq\left\lfloor n\left(1-x\right)\right\rfloor ,Y_{n}>\left\lfloor n\left(1+y\right)\right\rfloor \right)\\
 & = & \nu\left(Z_{\left\lfloor n\left(1+y\right)\right\rfloor }\leq\left\lfloor n\left(1-x\right)\right\rfloor \right)\\
 & = & \int_{K_{1}}\sum_{k=0}^{\left\lfloor n\left(1-x\right)\right\rfloor }1_{K_{1}\cap\left\{ \varphi>\left\lfloor n\left(1+y\right)\right\rfloor -k\right\} }\cdot\hat{T}^{k}\left(f\right)\; d\mu.\end{eqnarray*}
For $\delta\in\left(0,1-x\right)$ and $\varepsilon\in\left(0,1\right)$
fixed but arbitrary we divide the above sum into two parts as follows.\[
\nu\left(\frac{n-Z_{n}}{n}>x,\frac{Y_{n}-n}{n}>y\right)=\sum_{k=0}^{\left\lfloor n\delta\right\rfloor -1}\cdots+\sum_{k=\left\lfloor n\delta\right\rfloor }^{\left\lfloor n\left(1-x\right)\right\rfloor }\cdots=:I\left(n\right)+J\left(n\right).\]
By monotonicity of $\left(1_{K_{1}\cap\left\{ \varphi>n\right\} }\right)$
we have\begin{eqnarray*}
I\left(n\right) & \leq & \int_{K_{1}}1_{K_{1}\cap\left\{ \varphi>\left\lfloor n\left(1+y\right)\right\rfloor -\left\lfloor n\delta\right\rfloor +1\right\} }\cdot\sum_{k=0}^{\left\lfloor n\delta\right\rfloor -1}\hat{T}^{k}\left(f\right)\; d\mu.\end{eqnarray*}
Using (UA), (\ref{eq:AsymptAaronson}), and the fact that $K_{1}$
is uniform for $f$, we obtain for sufficiently large $n$ \begin{eqnarray*}
I\left(n\right) & \leq & \left(1+\varepsilon\right)^{2}\frac{\left\lfloor n\delta\right\rfloor -1}{\left\lfloor n\left(1+y\right)\right\rfloor -\left\lfloor n\delta\right\rfloor +1}\cdot\frac{1}{\log(\left\lfloor n\delta\right\rfloor -1)}\\
 & \sim & \left(1+\varepsilon\right)^{2}\frac{\delta}{y+1-\delta}\frac{1}{\log(n)}\quad\textrm{as}\; n\to\infty.\end{eqnarray*}
Thus, \[
\limsup_{n\to\infty}\log(n)\cdot I\left(n\right)\leq\left(1+\varepsilon\right)^{2}\frac{\delta}{y+1-\delta}.\]
Letting $\delta\to0$, we find\begin{equation}
I\left(n\right)=o\left(\frac{1}{\log(n)}\right),\quad\textrm{as}\; n\to\infty.\label{I(n)}\end{equation}
For the second part of the sum we have to show that\begin{equation}
J\left(n\right)\sim\frac{1}{\log\left(n\right)}\cdot\log\left(\frac{1+y}{x+y}\right)\qquad\textrm{as}\: n\to\infty.\label{J(n)}\end{equation}
Using the fact that $K_{1}$ is uniformly returning for $f$, a similarly
argument as in \cite{KesseboehmerSlassi:07}, Lemma 3.2, shows that
for all $n$ sufficiently large and $k\in\left[\left\lfloor n\delta\right\rfloor ,\left\lfloor n\left(1-x\right)\right\rfloor \right]$
we have uniformly on $K_{1}$\begin{equation}
\left(1-\varepsilon\right)\frac{1}{\log\left(n\right)}\leq\hat{T}^{k}\left(f\right)\leq\left(1+\varepsilon\right)^{2}\frac{1}{\log\left(n\right)}.\label{eq:(right W_n)}\end{equation}
Hence, using the right-hand side of (\ref{eq:(right W_n)}) and (\ref{Mu1}),
we obtain for $n$ sufficiently large\begin{eqnarray*}
J\left(n\right) & \leq & \frac{\left(1+\varepsilon\right)^{2}}{\log\left(n\right)}\cdot\sum_{k=\left\lfloor n\left(1+y\right)\right\rfloor -\left\lfloor n\left(1-x\right)\right\rfloor }^{\left\lfloor n\left(1+y\right)\right\rfloor -\left\lfloor n\delta\right\rfloor }\mu\left(K_{1}\cap\left\{ \varphi>k\right\} \right)\\
 & \sim & \left(1+\varepsilon\right)^{2}\frac{1}{\log\left(n\right)}\cdot\sum_{k=\left\lfloor n\left(1+y\right)\right\rfloor -\left\lfloor n\left(1-x\right)\right\rfloor }^{\left\lfloor n\left(1+y\right)\right\rfloor -\left\lfloor n\delta\right\rfloor }\frac{1}{k}\\
 & \sim & \left(1+\varepsilon\right)^{2}\frac{1}{\log\left(n\right)}\cdot\log\left(\frac{1+y-\delta}{x+y}\right)\quad\textrm{as}\; n\to\infty.\end{eqnarray*}
This implies\[
\limsup_{n\to\infty}\log\left(n\right)\cdot J\left(n\right)\leq\left(1+\varepsilon\right)^{3}\log\left(\frac{1+y-\delta}{x+y}\right).\]
Similarly, using the left-hand side of (\ref{eq:(right W_n)}), we
get \[
\liminf_{n\to\infty}\log\left(n\right)\cdot J\left(n\right)\geq\left(1+\varepsilon\right)^{2}\log\left(\frac{1+y-\delta}{x+y}\right).\]
Since $\varepsilon$ and $\delta$ were arbitrary, (\ref{J(n)}) holds.
Combining (\ref{I(n)}) and (\ref{J(n)}) proves our claim in (\ref{eqq1}).

\textbf{ad (2)} First, let $x\in\left(0,1\right)$. We have\begin{eqnarray*}
\nu\left(\frac{V_{n}}{n}>x\right) & = & \sum_{k=0}^{n}\nu\left(Y_{n}>k+\left\lfloor nx\right\rfloor ,Z_{n}=k\right)\\
 & = & \sum_{k=0}^{n-\left\lfloor nx\right\rfloor -1}\int_{K_{1}}1_{K_{1}\cap\left\{ \varphi>n-k\right\} }\cdot\hat{T}^{k}\left(f\right)\; d\mu\\
 &  & +\sum_{k=n-\left\lfloor nx\right\rfloor }^{n}\int_{K_{1}}1_{K_{1}\cap\left\{ \varphi>\left\lfloor nx\right\rfloor \right\} }\cdot\hat{T}^{k}\left(f\right)\; d\mu\\
 & =: & I\left(n\right)+J\left(n\right).\end{eqnarray*}
Let $\delta\in\left(0,1-x\right)$ and $\varepsilon\in\left(0,1\right)$
be fixed but arbitrary. First, we prove that\begin{equation}
J\left(n\right)\sim\frac{1}{\log(n)},\quad\textrm{as}\; n\to\infty.\label{(J(1))}\end{equation}
In fact, we have for sufficiently large $n$\begin{eqnarray*}
J\left(n\right) & \leq & \left(1+\varepsilon\right)^{2}\frac{\left\lfloor nx\right\rfloor -1}{\log(n)}\cdot\mu\left(K_{1}\cap\left\{ \varphi>\left\lfloor nx\right\rfloor \right\} \right)\\
 & \sim & \left(1+\varepsilon\right)^{2}\frac{1}{\log(n)}\quad\textrm{as}\; n\to\infty.\end{eqnarray*}
Similarly we get\begin{eqnarray*}
J\left(n\right) & \geq & \left(1-\varepsilon\right)\frac{\left\lfloor nx\right\rfloor -1}{\log(n)}\cdot\mu\left(K_{1}\cap\left\{ \varphi>\left\lfloor nx\right\rfloor \right\} \right)\\
 & \sim & \left(1-\varepsilon\right)\frac{1}{\log(n)}\quad\textrm{as}\; n\to\infty.\end{eqnarray*}
Combining both inequalities the asymptotic in (\ref{(J(1))}) follows. 

Now we have to prove that \begin{equation}
I\left(n\right)\sim-\log\left(x\right)\cdot\frac{1}{\log(n)}\qquad\textrm{as}\; n\to\infty.\label{II(n)}\end{equation}
We divide $I\left(n\right)$ again into two parts as follows\[
I\left(n\right)=\sum_{k=0}^{\left\lfloor n\delta\right\rfloor -1}\cdots+\sum_{k=\left\lfloor n\delta\right\rfloor }^{n-\left\lfloor nx\right\rfloor -1}\cdots=:I_{1}\left(n\right)+I_{2}\left(n\right).\]
Using the monotonicity of $\left(\1_{K_{1}\cap\left\{ \varphi>n\right\} }\right)$,
the fact that $K_{1}$ is uniformly for $f$, and (\ref{eq:AsymptAaronson})
we obtain for $n$ sufficiently large, \begin{eqnarray*}
I_{1}\left(n\right) & \leq & \int_{K_{1}}\1_{K_{1}\cap\left\{ \varphi>n-\left\lfloor n\delta\right\rfloor +1\right\} }\cdot\sum_{k=0}^{\left\lfloor n\delta\right\rfloor -1}\hat{T}^{k}\left(f\right)\; d\mu\\
 & \leq & \left(1+\varepsilon\right)^{2}\frac{\left\lfloor n\delta\right\rfloor -1}{n-\left\lfloor n\delta\right\rfloor +1}\cdot\frac{1}{\log(\left\lfloor n\delta\right\rfloor -1)}\\
 & \sim & \left(1+\varepsilon\right)^{2}\frac{\delta}{1-\delta}\frac{1}{\log(n)}\quad\textrm{as}\; n\to\infty.\end{eqnarray*}
Consequently,\begin{equation}
I_{1}\left(n\right)=o\left(\frac{1}{\log(n)}\right),\quad\textrm{as}\; n\to\infty.\label{eq:I_1(1)}\end{equation}
Now using the fact that $K_{1}$ is uniformly returning for $f$ we
have, for $n$ sufficiently large, \begin{eqnarray*}
I_{2}\left(n\right) & \leq & \frac{\left(1+\varepsilon\right)^{2}}{\log(n)}\cdot\sum_{k=\left\lfloor nx\right\rfloor +1}^{n-\left\lfloor n\delta\right\rfloor }\mu\left(K_{1}\cap\left\{ \varphi>k\right\} \right)\\
 & \sim & \left(1+\varepsilon\right)^{2}\frac{1}{\log(n)}\cdot\log\left(\frac{1-\delta}{x}\right)\quad\textrm{as}\; n\to\infty.\end{eqnarray*}
This implies \[
\limsup_{n\to\infty}\log(n)\cdot I_{2}\left(n\right)\leq\left(1+\varepsilon\right)^{3}\log\left(\frac{1-\delta}{x}\right).\]
Similarly, we get \[
\liminf_{n\to\infty}\log(n)\cdot I_{2}\left(n\right)\geq\left(1+\varepsilon\right)^{2}\log\left(\frac{1-\delta}{x}\right).\]
Since $\varepsilon$ and $\delta$ were arbitrary, we have\begin{equation}
I_{2}\left(n\right)\sim-\log\left(x\right)\cdot\frac{1}{\log(n)},\quad\textrm{as}\: n\to\infty.\label{eq:I_2(1)}\end{equation}
The asymptotics (\ref{eq:I_1(1)}) and (\ref{eq:I_2(1)}) prove (\ref{II(n)}).
Combining (\ref{(J(1))}) and (\ref{II(n)}) we obtain for $x\in\left(0,1\right)$\begin{equation}
\nu\left(\frac{V_{n}}{n}>x\right)\sim\frac{1-\log(x)}{\log(n)}\quad\textrm{as}\; n\to\infty.\label{Vn1}\end{equation}

Now we consider the case $x\geq1$. Since\begin{eqnarray*}
\nu\left(\frac{V_{n}}{n}>x\right) & = & \sum_{k=0}^{n}\int_{K_{1}}\1_{K_{1}\cap\left\{ \varphi>\left\lfloor nx\right\rfloor \right\} }\cdot\hat{T}^{k}\left(f\right)\; d\mu\end{eqnarray*}
we have, for $n$ sufficiently large, \begin{eqnarray*}
\nu\left(\frac{V_{n}}{n}>x\right) & \leq & \mu\left(K_{1}\cap\left\{ \varphi>\left\lfloor nx\right\rfloor \right\} \right)\frac{n}{\log(n)}\left(1+\varepsilon\right)\\
 & \sim & \frac{1}{x\log(n)}\left(1+\varepsilon\right).\end{eqnarray*}
Similarly, we obtain the reverse inequality proving the statement
that for $x\geq1$\begin{equation}
\nu\left(\frac{V_{n}}{n}>x\right)\sim\frac{1}{x\log(n)}\quad\textrm{as}\; n\to\infty.\label{Vn2}\end{equation}
Combining (\ref{Vn1}) and (\ref{Vn2}) the statement (\ref{eqq2})
follows from (\ref{teta2}).

\textbf{ad (3)} We first observe from (\ref{teta1}) and (\ref{teta2})
that on $\I$ we have \[
\frac{\kappa_{\theta_{n}+1}}{\sum_{k=1}^{\theta_{n}+1}\kappa_{k}}=\frac{V_{n-1}}{Y_{n-1}}\quad\textrm{for\; all}\; n\geq1.\]
For $x\in\left(0,1\right)$we have\begin{eqnarray*}
\nu\left(\frac{V_{n}}{n}>x\right) & = & \sum_{k=0}^{n}\nu\left(Y_{n}>k+\left\lfloor \frac{x}{1+x}k\right\rfloor ,Z_{n}=k\right)\\
 & = & \sum_{k=0}^{\left\lfloor n\left(1-x\right)\right\rfloor }\int_{K_{1}}\1_{K_{1}\cap\left\{ \varphi>n-k\right\} }\cdot\hat{T}^{k}\left(f\right)\; d\mu\\
 &  & +\sum_{k=\left\lfloor n\left(1-x\right)\right\rfloor +1}^{n}\int_{K_{1}}\1_{K_{1}\cap\left\{ \varphi>\left\lfloor \frac{x}{1+x}k\right\rfloor \right\} }\cdot\hat{T}^{k}\left(f\right)\; d\mu.\end{eqnarray*}
Then a similarly argument as in \textbf{ad (2)} gives the statement
in (\ref{eqq3}).

\textbf{ad (4)} Clearly, (\ref{eqq4}) is a direct consequence from
(\ref{eqq3}) using the fact that on $\left\{ \theta_{n}>0\right\} $
we have\[
\frac{\kappa_{\theta_{n}+1}}{\sum_{k=1}^{\theta_{n}}\kappa_{k}}=\frac{\kappa_{\theta_{n}+1}/\sum_{k=1}^{\theta_{n}+1}\kappa_{k}}{1-\kappa_{\theta_{n}+1}/\sum_{k=1}^{\theta_{n}+1}\kappa_{k}}.\]

\end{proof}

\end{document}